\documentclass{elsart}
\usepackage{amssymb,amsmath}
\usepackage{latexsym}

\begin{document}

\title{Applications of tangent transformations to the linearization problem of fourth-order
ordinary differential equations}

\author[NU]{S. Suksern},
\ead{supapornsu@nu.ac.th}
\author[SUT]{S.V. Meleshko}

\address[NU]{Department of Mathematics, Faculty of Science, Naresuan University,
Phitsanulok, 65000, Thailand}
\address[SUT]{School of Mathematics, Suranaree University of Technology,
Nakhon Ratchasima, 30000, Thailand}

\begin{abstract}
Linearization problem of ordinary differential equations by a
new set of tangent transformations is considered in the paper. This
set of transformations allows one to extend the set of transformations
applied for the linearization problem. Criteria for fourth-order
ordinary differential equations to be linearizable are obtained in
a particular case; a linearization test and procedure for
obtaining the linearizing transformations are provided in explicit
forms.
\end{abstract}


\section{Introduction}

\label{Intro}

The problem of linearization of a nonlinear ordinary differential
equation has attracted a lot of attention of scientists. The first
linearization problem for ordinary differential equations was solved
by S.Lie \cite{bk:Lie[1883]}. He found the general form of all ordinary
differential equations of second-order which can be reduced to a linear
equation by changing the independent and dependent variables. He showed
that any linearizable second-order ordinary differential equation should be at most cubic
in the first-order derivative, and provided a linearization test in
terms of its coefficients.

It should be noted that the linearization problem is a particular case of an
equivalence problem, where one has to map one equation into another.
A different approach for tackling an equivalence problem of second-order
ordinary differential equations was developed by E.Cartan \cite{bk:Cartan[1924]}.
The idea of his approach was to associate with every differential
equation a uniquely defined geometric structure of a certain form.

These two approaches were also applied to third- and fourth-order ordinary
differential equations. The linearization problem of a third-order
ordinary differential equation with respect to point transformations
was studied in \cite{bk:Grebot[1997]}. Complete criteria for linearization
of third-order were obtained in \cite{bk:IbragimovMeleshko[2004],bk:IbragimovMeleshko[2005]}
and fourth-order in \cite{bk:IbragimovMeleshkoSuksern[2008]}.

S.Lie also noted that all second-order ordinary differential equations
can be transformed to each other by means of contact transformations.
Because this is not so for ordinary differential equations of higher
order, contact transformations became interesting for applying them to
a linearization problem. For third-order ordinary differential equations
contact transformations were studied in \cite{bk:Chern[1940],bk:BocharovSokolovSvinolupov[1993],bk:GusyatnikovaYumaguzhin[1999],bk:DoubrovKomrakovMorimoto[1999],bk:Doubrov[2001]}.
The complete solution of the linearization problem was given in \cite{bk:NeutPetitot[2002]}
and in explicit form the criteria for linearization is presented in
\cite{bk:IbragimovMeleshko[2005]}. For equations of fourth-order
the linearization problem via contact transformations was considered
in \cite{bk:DridiNeut[2005]} and a complete solution is given in
\cite{bk:SuksernIbragimovMeleshko[2009]}.

It is known that point and contact transformations are exceptional
types of tangent transformations
\begin{equation}
t=\varphi(x,y,y^{\prime},...,y^{(s)}),\ \ u=\psi(x,y,y^{\prime},...,y^{(s)}).\label{eq:tangent_tr}
\end{equation}
 Here for the purpose of the present paper a single independent variable
$x$ and a single dependent variable $y$ are used. Notice that
the requirement to satisfy the tangent conditions defines the transformations
of the derivatives through the functions $\varphi(x,y,y^{\prime},...,y^{(s)})$
and $\psi(x,y,y^{\prime},...,y^{(s)})$:
\[
u^{(j)}=\frac{Du^{(j-1)}}{D\varphi},\;(j\geq1),
\]
 where $u^{(0)}=\psi$, and $D$ is the operator of the total derivative
with respect to $x$:
\[
D=\partial_{x}+\sum\limits _{j=0}y^{(j+1)}\partial_{y^{(j)}},
\]
 with $y^{(0)}=y$. Remind also that if the function $y=y_{o}(x)$
is given, then the transformed function $u_{o}(t)$ is defined by
the following two steps. On the first step one has to solve with respect
to $x$ the equation
\[
t=\varphi(x,y_{o}(x),y_{o}^{\prime}(x),...,y_{o}^{(s)}(x)).
\]
 Let say $x=g(t)$ is a solution of this equation. The transformed
function is determined by the formula
\[
u_{o}(t)=\varphi(g(t),y_{o}(g(t)),y_{o}^{\prime}(g(t)),...,y_{o}^{(s)}(g(t))).
\]
 The inverse transformation is more complicated. If one has the function
$u_{o}(t)$, then for finding the function $y_{o}(x)$ one has to
solve the ordinary differential equation
\[
u_{o}\left(\varphi(x,y_{o}(x),y_{o}^{\prime}(x),...,y_{o}^{(s)}(x))\right)=\psi(x,y_{o}(x),y_{o}^{\prime}(x),...,y_{o}^{(s)}(x)).
\]

The set of point transformations corresponds to $s=0$. In this case
for the forward and inverse transformations one needs to overcome
similar difficulties. For the set of contact transformations $s=1$,
and it is also required that $u^{\prime}$ do not depend on second-order
derivatives $y^{\prime\prime}$. In this case the forward transformation
is simpler than the inverse transformation since for the forward transformation
one has just to use the differentiation, whereas for the inverse transformation
one needs to solve the first-order ordinary differential equation
\begin{equation}
u_{o}\left(\varphi(x,y_{o}(x),y_{o}^{\prime}(x))\right)=\psi(x,y_{o}(x),y_{o}^{\prime}(x)).\label{eq:inverse_cont}
\end{equation}

The sets of point and contact transformations are exceptional because
the transformed derivative $u^{(j)}$ depends on derivatives of order
not higher than $j$. Moreover, B\"acklund \cite{bk:Backlund[1874]}
proved that the tangent transformations closed in a finite-dimensional
space of the independent, dependent variables and derivatives of finite
order (say order $s$) are comprised by point and contact transformations%
\footnote{This statement is called the B\"acklund non-existence theorem \cite{bk:HandbookLie_v3,bk:Ibragimov[1983],bk:Ibragimov[1999]}.%
}. It has to be mentioned here the following.

First, the B\"acklund theorem was proven with the assumption that there
are no restrictions on derivatives except satisfying the tangent conditions,
whereas it is not so for a set of solutions of differential equations.
This circumstance gives possibilities for the existence of tangent
transformations of finite order which act on the set of solutions%
\footnote{Details of this discussion can be found in \cite{bk:Meleshko[2005]}.}.

Second, for a linearization problem a change of order of an equation
is not crucial. There are attempts to tackle the linearization problem
by increasing order of an equation \cite{bk:FerapontovSvirshchevskii[2007],bk:AndriopoulosLeach[2007-1],bk:IbragimovMeleshko[2007-1]},
reducing order \cite{bk:Meleshko[2006]lin} or their combinations
\cite{bk:Abraham-Shrauner[1993]a}.

Third, difficulties in applications of the inverse transformation
of contact transformations and tangent transformations with $s=1$
are similar. This allows one to assume that tangent transformations
which are not contact could be also useful for a linearization problem.

Notice also that through involving derivatives in a
transformation, the notion of a Lie group of point transformations
is generalized. The limitations dictated by the B\"acklund theorem
can be overcome by
considering the admitted transformations%
\footnote{Examples of such transformations one can find in \cite{bk:Meleshko[2005]}.%
} or extending the space of the derivatives involved in the transformations
up to infinity. The second alternative leads to Lie-B\"acklund transformations
\cite{bk:Ibragimov[1983]}.

The present paper deals with tangent transformations (\ref{eq:tangent_tr})
of first-order ($s=1$) which are not point or contact transformations.
It is worth to notice that the classical reduction of order and the
Ricatti substitution can be also considered as a particular case of
tangent transformations of the form (\ref{eq:tangent_tr}):
\[
t=\varphi(x,y,p),\quad u=\psi(x,y,p),
\]
 where $p=y^{\prime}$. Indeed, for the classical reduction of order
the transformation is
\[
t=y,\;\; u=p,
\]
 and for the Ricatti substitution
\[
t=x,\;\; u=p/y.
\]
 For a linearization problem these transformations were applied in
\cite{bk:Meleshko[2006]lin} and in \cite{bk:IbragimovMeleshko[2007-1]}.

The main goal of the present paper is to demonstrate possibilities of applications
tangent transformations for a linearization problem.
In the paper tangent transformations are applied for linearizing
fourth-order ordinary differential equations. Complete study of fiber
preserving transformations ($\varphi_{p}=\varphi_{y}=0$) mapping equations to the trivial
third-order ordinary differential equation $y^{(3)}=0$ is given in the paper.

The manuscript is organized as follows. In Section 2, the necessary conditions of
linearization of a fourth-order ordinary differential
equation in the general case of tangent transformations with $s=1$ are presented.
It is shown that linearizing equations are separated in two classes of equations.
Further study is devoted to one of these classes. In particular, in Section 3,
sufficient conditions for linearizing fiber preserving tangent transformations are obtained.
In Section 4, we state the theorems that yield criteria for a fourth-order ordinary differential equation to
be linearizable via fiber preserving tangent transformations. Relations between coefficients of a linearizable
equation and tangent transformations
reducing this equation into the equation $y^{(3)}=0$ are presented in this section. These relations are necessary for
formulating the linearization theorems which are given in Section 4.
For the sake of simplicity of reading, cumbersome formulae are moved into Appendix.
Examples which illustrate the procedure of using
the linearization theorems are presented in Section 5.

\section{Necessary conditions of linearization}

We begin with obtaining necessary conditions for linearizable equations.
Recall that according to the Laguerre theorem a linear third-order
ordinary differential equation has the form
\begin{equation}
u^{\prime\prime\prime}+\alpha(t)u=0.\label{eq:03}
\end{equation}
 The necessary form of a linearizable fourth-order ordinary differential
equation
\begin{equation}
y^{\left(4\right)}=f(x,y,y^{\prime},y^{\prime\prime},y^{\prime\prime\prime})\label{eq:01}
\end{equation}
 is obtained by the substitution of transformed derivatives into the
Laguerre form (\ref{eq:03}). Using the tangent conditions, the transformation
\begin{equation}
t=\varphi(x,y,p),\quad u=\psi(x,y,p),\label{eq:02}
\end{equation}
 defines the change of derivatives
\begin{equation}
u^{\prime}=\frac{D\psi}{D\varphi},\; u^{\prime\prime}=\frac{Du^{\prime}}{D\varphi},\; u^{\prime\prime\prime}=\frac{Du^{\prime\prime}}{D\varphi}.\label{eq:04}
\end{equation}
 Substituting the resulting expressions of the derivatives into linear
equation (\ref{eq:03}), one arrives at the following equation%
\footnote{Since the study requires a huge amount of analytical calculations,
it is necessary to use a computer for these calculations. A brief
review of computer systems of symbolic manipulations can be found,
for example, in \cite{bk:Davenport[1994]}. In our calculations the
system REDUCE \cite{bk:Hearn} was used.%
}
\begin{equation}
(\varphi_{p}(\psi_{x}+p\psi_{y})-\psi_{p}(\varphi_{x}+p\varphi_{y}))y^{(4)}+H(x,y,p,y^{\prime\prime},y^{\prime\prime\prime})=0,\label{eq:07}
\end{equation}
 with some function $H(x,y,y^{\prime},y^{\prime\prime},y^{\prime\prime\prime})$.
Since for
\begin{equation}
\psi_{p}(\varphi_{x}+p\varphi_{y})=\varphi_{p}(\psi_{x}+p\psi_{y})\label{eq:06}
\end{equation}
 the tangent transformation is a contact transformation%
\footnote{Complete study of this case is given in \cite{bk:NeutPetitot[2002]}
and \cite{bk:IbragimovMeleshko[2005]}.%
}, the case $\varphi_{p}(\psi_{x}+p\psi_{y})-\psi_{p}(\varphi_{x}+p\varphi_{y})\neq0$
is considered in the paper. The study of the function $H(x,y,y^{\prime},y^{\prime\prime},y^{\prime\prime\prime})$
in (\ref{eq:07}) leads to the property that the transformations (\ref{eq:02})
with $\varphi_{p}=0$ and $\varphi_{p}\not=0,$ provide two distinctly
different candidates for the linearizable equations:\\
 (a) if $\varphi_{p}=0$, then the linearizable equation has the form
\begin{equation}
y^{(4)}+(A_{1}y^{\prime\prime}+A_{0})y^{\prime\prime\prime}+B_{3}y^{\prime\prime}\,^{3}+B_{2}y^{\prime\prime}\,^{2}+B_{1}y^{\prime\prime}+B_{0}=0,\label{eq:08}
\end{equation}
 (b) if $\varphi_{p}\not=0$, one arrives to the equation
\begin{equation}
\begin{array}{c}
{\displaystyle y^{\left(4\right)}+\frac{1}{{y^{\prime\prime}+r}}\Big[-3{y^{\prime\prime\prime}}^{2}+\left({C_{2}y^{\prime\prime}{}^{2}+C_{1}y^{\prime\prime}+C_{0}}\right)y^{\prime\prime\prime}}\\
+D_{5}y^{\prime\prime}{}^{5}+D_{4}y^{\prime\prime}{}^{4}+D_{3}y^{\prime\prime}{}^{3}+D_{2}y^{\prime\prime}{}^{2}+D_{1}y^{\prime\prime}+D_{0}\Big]=0.
\end{array}\label{eq:09}
\end{equation}
 Here $A_{i}=A_{i}(x,y,p)$, $B_{i}=B_{i}(x,y,p)$, $C_{i}=C_{i}(x,y,p)$
and $D_{i}=D_{i}(x,y,p)$ are some functions of $x,y,p$. Expressions
of the coefficients $A_{i}=A_{i}(x,y,p)$, $B_{j}=B_{j}(x,y,p)$,
($i=0,1;\ j=0,1,2,3$) in the case $\alpha(t)=0$ are presented in
Appendix.

Thus, it is shown that every fourth-order ordinary differential equation
linearizable by a tangent transformation of the form (\ref{eq:02})
belongs either to the class of (\ref{eq:08}) or to the class of (\ref{eq:09}).



\section{Sufficient conditions of linearization for the case (\ref{eq:08})}

For obtaining sufficient conditions, one has to solve the compatibility
problem considering the representations of the coefficients $A_{i}(x,y,p)$,
$B_{j}(x,y,p)$ through the unknown functions $\varphi$, $\psi$
and $\alpha$ with the given coefficients $A_{i}(x,y,p)$, $B_{j}(x,y,p)$.
For example, in the case $\alpha=0$, these equations are (\ref{ap:01})-(\ref{ap:06}).

For simplicity of calculations the function $k=\varphi_{x}+p\varphi_{y}\neq0$
is introduced. From this definition one gets the derivative $\varphi_{x}=k-p\varphi_{y}$.
Notice that the function $k$ is a linear function with respect to
$p$ with $k_{p}=\varphi_{y}$.

From equations (\ref{ap:01}), (\ref{ap:02}) and (\ref{ap:04}) one
can find the derivatives $\psi_{pp}$, $\psi_{px}$ and $\psi_{py}$,
respectively. Substituting them into equation (\ref{ap:03}), it becomes
\begin{equation}
\varphi_{y}=k\lambda,\label{eq:may21.1}
\end{equation}
 where
\begin{equation}
5\lambda^{2}+A_{1}\lambda-(3A_{1p}+A_{1}^{2}-9B_{3})=0.\label{eq:may25.1a}
\end{equation}
 The discriminant of this quadratic equation has to be not negative
\begin{equation}
\Delta\equiv20A_{1p}+7A_{1}^{2}-60B_{3}\geq0.\label{eq:discriminant}
\end{equation}
 Defining the function $q(x,y,p)$ such that
\begin{equation}
q^{2}=3\Delta\label{eq:may25.1}
\end{equation}
 one finds
\[
\lambda=-(A_{1}+q)/10.
\]
 Differentiating equation (\ref{eq:may21.1}) with respect to $p$,
one finds that
\begin{equation}
\lambda_{p}+\lambda^{2}=0.\label{eq:jun7.1}
\end{equation}
 Comparing the mixed derivatives $(\varphi_{y})_{x}=(\varphi_{x})_{y}$,
one gets
\begin{equation}
\lambda k_{x}+(p\lambda-1)k_{y}+(\lambda_{x}+p\lambda_{y})k=0.\label{ad1}
\end{equation}

Because of the representation of equations (\ref{eq:jun7.1}) and
(\ref{ad1}), further analysis of the compatibility depends on value
of $\lambda$. In the present paper a complete solution of the case
where $\lambda=0$ is given. In this case $\varphi=\varphi(x)$. Transformations
of the class where $\varphi=\varphi(x)$ are called fiber preserving
transformations and they are actively applied to a linearization problem
\cite{bk:Grebot[1997],bk:Chern[1940],bk:BocharovSokolovSvinolupov[1993]}.

Assuming that $\lambda=0$, from equation (\ref{ap:05}) one finds the derivative
$\psi_{xy}$. From the relation $(\psi_{py})_{x}-(\psi_{px})_{y}=0$,
one obtains the derivative $\psi_{yy}$. The equation $(\psi_{px})_{p}-(\psi_{pp})_{x}=0$
provides the relation
\begin{equation}
\mu_{2}\psi_{p}-A_{1}\psi_{y}=0,\label{newmm2}
\end{equation}
 where
\[
\mu_{1}=(3A_{0p}+A_{0}A_{1}-3B_{2})/3,
\]
\[
\mu_{2}=3A_{1x}+3A_{1y}p+A_{0}A_{1}-3B_{2}+9\mu_{1}.
\]

Further analysis of the compatibility depends on value of $A_{1}$:
it is separated into two cases $A_{1}=0$ and $A_{1}\neq0$. Notice
also that equation (\ref{ap:01}) becomes
\begin{equation}
\psi_{pp}=A_{1}\psi_{p}/3.\label{eq:dec08.2012.1}
\end{equation}


\subsection{Case $A_{1}\neq0$}

Using equation (\ref{newmm2}), one can exclude the derivative $\psi_{y}$
from other equations. Comparing $(\psi_{y})_{y}$ with the found derivative $\psi_{yy}$, and
composing the relations $(\psi_{y})_{p}=\psi_{py}$
and $(\psi_{py})_{p}=(\psi_{pp})_{y}$, one has the conditions
\begin{equation}
\mu_{2y}=(A_{0y}A_{1}^{2}+A_{1y}\mu_{2}-3\mu_{1x}A_{1}^{2}-3\mu_{1y}A_{1}^{2}p)/A_{1},\label{suff:01}
\end{equation}
\begin{equation}
\mu_{1p}=A_{1y}/3,\label{suff:02}
\end{equation}
\begin{equation}
\begin{array}{c}
B_{2p}=(6A_{1y}A_{1}+9B_{3x}A_{1}+9B_{3y}A_{1}p+3A_{0}A_{1}B_{3}\\
-A_{1}^{2}B_{2}+6A_{1}^{2}\mu_{1}-3B_{3}\mu_{2})/(3A_{1}).
\end{array}\label{suff:03}
\end{equation}
 The equation $(\psi_{px})_{y}=(\psi_{y})_{xp}$ can be solved with
respect to $\varphi_{xxx}$:
\begin{equation}
\varphi_{xxx}=(3\varphi_{xx}^{2}+2\varphi_{x}^{2}\mu_{3})/(2\varphi_{x}),\label{tran:02}
\end{equation}
 where
\[
\begin{array}{rl}
\mu_{3}= & (-9A_{0y}A_{1}+3B_{1p}A_{1}-3B_{2x}A_{1}-3B_{2y}A_{1}p+9\mu_{1x}A_{1}+9\mu_{1y}A_{1}p\\
 & -A_{0}A_{1}B_{2}-6A_{0}A_{1}\mu_{1}+A_{1}^{2}B_{1}+B_{2}\mu_{2}-3\mu_{1}\mu_{2})/(2A_{1}^{2}).
\end{array}
\]
 Since $\varphi=\varphi(x)$, one has
\begin{equation}
\mu_{3y}=0,\ \ \ \mu_{3p}=0.\label{suff:04}
\end{equation}
 The relation $(\psi_{y})_{x}=\psi_{xy}$ gives the condition
\begin{equation}
\begin{array}{rl}
\mu_{2x}= & (-3A_{0x}A_{1}^{2}-9A_{0y}A_{1}^{2}p-6A_{1y}\mu_{2}p+18\mu_{1x}A_{1}^{2}p\\
 & +18\mu_{1y}A_{1}^{2}p^{2}-A_{0}^{2}A_{1}^{2}-3A_{0}A_{1}\mu_{2}+3A_{1}^{2}B_{1}-6A_{1}^{2}\mu_{3}\\
 & +6B_{2}\mu_{2}-18\mu_{1}\mu_{2}+4\mu_{2}^{2})/(6A_{1}).
\end{array}\label{suff:06}
\end{equation}
 Equation (\ref{ap:06}) can be solved with respect to
\begin{equation}
\psi_{xxx}=(-3\varphi_{xx}^{2}\psi_{x}+6\varphi_{xx}\varphi_{x}\psi_{xx}+\varphi_{x}^{2}\psi_{p}\mu_{4}+2\varphi_{x}^{2}\psi_{x}\mu_{3})/(2\varphi_{x}^{2}),\label{tran:03}
\end{equation}
 where
\begin{equation*}
    \begin{array}{rl}
    \mu_{4}= &(6 A_{0xy} A_{1} p^2 + 6 A_{0xx} A_{1} p + 7 A_{0x} A_{0} A_{1} p - 6 A_{0x} A_{1} \mu_{1} p^2 - 4 A_{0x} \mu_{2} p \\
    &- 4 A_{0yy} A_{1} p^3+ 5 A_{0y} A_{0} A_{1} p^2 - 14 A_{0y} A_{1} \mu_{1} p^3 - 2A_{0y} \mu_{2} p^2 - 6 B_{1x} A_{1} p \\
    &+ 12 \mu_{1xy} A_{1} p^3 - 6\mu_{1x} A_{0} A_{1} p^2 + 24 \mu_{1x} A_{1} \mu_{1} p^3 + 12\mu_{1yy} A_{1} p^4 \\
    &- 6 \mu_{1y} A_{0} A_{1} p^3 + 24 \mu_{1y}A_{1} \mu_{1} p^4 - 4 \mu_{1y} \mu_{2} p^3 + 12 \mu_{3x} A_{1} p \\
    &+ A_{0}^3 A_{1} p - 2 A_{0}^2 A_{1} \mu_{1} p^2 - A_{0}^2 \mu_{2} p- 3 A_{0} A_{1} B_{1} p + 6 A_{0} A_{1} \mu_{3} p \\
    &+ 2 A_{0}\mu_{1} \mu_{2} p^2 + 4 A_{1} B_{0} + 6 A_{1} B_{1} \mu_{1} p^2 -12 A_{1} \mu_{1} \mu_{3} p^2 + 2 B_{1} \mu_{2} p \\
    &- 4 \mu_{1}^2\mu_{2} p^3 - 12 \mu_{2} \mu_{3} p)/(2 A_{1}).
    \end{array}
\end{equation*}
 Comparing the mixed derivative $(\psi_{xxx})_{p}=(\psi_{px})_{xx}$
and $(\psi_{xxx})_{y}=(\psi_{y})_{xxx}$, one obtains conditions (\ref{suff:07}) and (\ref{suff:08}).

The obtained system of equations for the functions $\varphi(x)$ and
$\psi(x,y,p)$ consisting of the equations
\begin{equation}\label{tran:04}
    \psi_{pp} = (\psi_{p} A_{1})/3,
\end{equation}
\begin{equation}\label{tran:05}
    \psi_{px} = (3 \varphi_{xx} \psi_{p} A_{1} + \varphi_{x} \psi_{p} (A_{0} A_{1} - 3 A_{1} \mu_{1} p - \mu_{2}))/(3 \varphi_{x} A_{1}),
\end{equation}
\begin{equation}\label{tran:01}
    \psi_{y} = (\psi_{p} \mu_{2})/A_{1},
\end{equation}
and satisfying the condition $(\psi _{py})_{x}=(\psi _{xy})_{p}$,
is involutive.

\subsection{Case $A_{1}=0$}

From equation (\ref{newmm2}), one obtains the condition
\begin{equation}
\mu_{2}=0,\label{suff:09}
\end{equation}
 which means that
\begin{equation}
\mu_{1}=B_{2}/3,\label{suff:11}
\end{equation}
 Substituting $A_{1}$ and $\mu_{1}$ into (\ref{eq:may25.1a}) and
the definition of $\mu_{1}$, one arrives at the conditions
\begin{equation}
B_{3}=0,\,\,\, A_{0p}=4B_{2}/3.\label{suff:10}
\end{equation}
 Equating the mixed derivatives $(\psi_{py})_{p}-(\psi_{pp})_{y}=0$,
one gets $B_{2}=B_{2}(x,y)$. Hence, one finds that
\begin{equation}
A_{0}=4pB_{2}/3+A_{00},\label{suff:12}
\end{equation}
 where $A_{00}(x,y)$ is a function of the integration. Equating the mixed
derivatives $(\psi_{px})_{y}=(\psi_{xy})_{p}$, and integrating it
with respect to $p$, one derives
\begin{equation}
B_{1}=(9A_{00y}p+6B_{2y}p^{2}+3A_{00}B_{2}p+2B_{2}^{2}p^{2})/3+B_{10},\label{suff:13}
\end{equation}
 where $B_{10}=B_{10}(x,y)$.

Since $\psi_{pp}=0$, then
\begin{equation}
\psi=\psi_{1}p+\psi_{0},\label{tran:06}
\end{equation}
 where $\psi_{0}(x,y)$ and $\psi_{1}(x,y)\neq0$ are some functions.
Substituting (\ref{tran:06}) into the relation for $\psi_{py},\psi_{px}$
and $\psi_{xy}$, one finds that
\begin{equation}
\psi_{1y}=B_{2}\psi_{1}/3,\label{tran:07}
\end{equation}
 and the derivatives $\psi_{0y}$ and $\varphi_{xxx}$, respectively.

Equation (\ref{ap:06}) becomes a polynomial of degree 3 with respect
to $p$ with coefficients not depending on $p$. One can split this
equation. Equating the coefficient with respect to $p$, and integrating
it, one obtains the relation:
\begin{equation}
\begin{array}{rl}
B_{0}= & (27A_{00yy}p^{3}+18A_{00y}B_{2}p^{3}+9B_{2yy}p^{4}+9B_{2y}A_{00}p^{3}+9B_{2y}B_{2}p^{4}\\
 & +3A_{00}B_{2}^{2}p^{3}+B_{2}^{3}p^{4})/27+B_{00}+B_{01}p+B_{02}p^{2},
\end{array}\label{suff:14}
\end{equation}
 where $B_{00}(x,y)$, $B_{01}(x,y)$ and $B_{02}(x,y)$ are some
functions. From the other three coefficients one can find $\psi_{0xxx}$,
$\psi_{1xxx}$ and
\begin{equation}
B_{10y}=(3B_{02}-B_{10}B_{2})/3.\label{suff:15}
\end{equation}

Comparing the mixed derivatives $(\psi_{yy})_{x}=(\psi_{xy})_{y}$,
one arrives at the equation
\begin{equation}
\mu_{5}(\varphi_{xx}\psi_{1}-\varphi_{x}\psi_{1x}+\varphi_{x}\mu_{6}\psi_{1})=0,\label{eq:10}
\end{equation}
 where
\[
\mu_{5}=3(-3A_{00y}+4B_{2x}),\,\,\,\,\mu_{5x}+\mu_{5}\mu_{6}=6B_{2xx}+3B_{2x}A_{00}-3B_{02}+B_{10}B_{2}.
\]

\subsubsection{Case $\mu_{5}\neq0$}

From equation (\ref{eq:10}), one finds
\begin{equation}
\varphi_{xx}=(\varphi_{x}\psi_{1x}-\varphi_{x}\mu_{6}\psi_{1})/\psi_{1}.\label{tran:08}
\end{equation}
 Differentiating this equation with respect to $y$ and $x$, one
gets
\begin{equation}
3\mu_{6y}-B_{2x}=0,\label{suff:16}
\end{equation}
\begin{equation}
\psi_{1xx}=(3\psi_{1x}^{2}-2\psi_{1x}\mu_{6}\psi_{1}+\psi_{1}^{2}(-3A_{00x}+8\mu_{6x}-3A_{00}\mu_{6}+B_{10}+7\mu_{6}^{2}))/(2\psi_{1}).\label{tran:09}
\end{equation}

Substituting $\psi_{1xx}$ into the found earlier expression for $\psi_{1xxx}$,
one obtains the condition
\begin{equation}
\begin{array}{rl}
\mu_{6xx}= & (3A_{00xx}-6A_{00x}A_{00}+21A_{00x}\mu_{6}+B_{10x}+15\mu_{6x}A_{00}-60\mu_{6x}\mu_{6}\\
 & -6A_{00}^{2}\mu_{6}+2A_{00}B_{10}+30A_{00}\mu_{6}^{2}-2B_{01}-4B_{10}\mu_{6}-40\mu_{6}^{3})/10.
\end{array}\label{suff:17}
\end{equation}
Comparing the mixed derivative $(\psi_{1xxx})_{y}=(\psi_{1y})_{xxx}$, one arrives at the condition
\begin{equation}
\begin{array}{rl}
\mu_{5xx} = &( - 36 A_{00x} \mu_{6y} + 6 A_{00x} \mu_{5} - 18 B_{01y} + 12 B_{02x} - 4 B_{10x} B_{2} + \mu_{5x} A_{00} - 8 \mu_{5x} \mu_{6}\\
&- 14 \mu_{6x} \mu_{5} - 9 \mu_{6y} A_{00}^2 + 24 \mu_{6y} B_{10} + 3 A_{00} B_{02} - A_{00} B_{10} B_{2} + 7 A_{00} \mu_{5} \mu_{6} \\
&- 2 B_{10} \mu_{5} - 18 \mu_{5} \mu_{6}^2)/2.
\end{array}\label{suff:26}
\end{equation}
 The equation $(\psi_{0xxx})_{y}=(\psi_{0y})_{xxx}$ yields the condition
\begin{equation}
\begin{array}{c}
-9A_{00xx}A_{00}-18A_{00xx}\mu_{6}+72A_{00x}^{2}-432A_{00x}\mu_{6x}-18A_{00x}A_{00}^{2}\\
+189A_{00x}A_{00}\mu_{6}-24A_{00x}B_{10}-432A_{00x}\mu_{6}^{2}+60B_{00y}-36B_{01x}\\
+18B_{10xx}+9B_{10x}A_{00}+36B_{10x}\mu_{6}+540\mu_{6x}^{2}+27\mu_{6x}A_{00}^{2}\\
-540\mu_{6x}A_{00}\mu_{6}+108\mu_{6x}B_{10}+1080\mu_{6x}\mu_{6}^{2}-18A_{00}^{3}\mu_{6}+6A_{00}^{2}B_{10}\\
+162A_{00}^{2}\mu_{6}^{2}-6A_{00}B_{01}-36A_{00}B_{10}\mu_{6}-540A_{00}\mu_{6}^{3}+20B_{00}B_{2}\\
-36B_{01}\mu_{6}+108B_{10}\mu_{6}^{2}+540\mu_{6}^{4}-6A_{00xxx}=0.
\end{array}\label{suff:19}
\end{equation}

The obtained system of equations for the functions $\varphi(x)$,
$\psi_{1}(x,y)$ and $\psi_{0}(x,y)$ consisting of the equations
\begin{equation}
\psi_{0y}=\psi_{1}(A_{00}-3\mu_{6}),\label{tran:10}
\end{equation}
\begin{equation}
\begin{array}{rl}
\psi_{0xxx}= & (6\psi_{0xx}\psi_{1x}\psi_{1}-6\psi_{0xx}\mu_{6}\psi_{1}^{2}-3\psi_{0x}\psi_{1x}^{2}+6\psi_{0x}\psi_{1x}\mu_{6}\psi_{1}\\
 & +\psi_{0x}\psi_{1}^{2}(-3A_{00x}+6\mu_{6x}-3A_{00}\mu_{6}+B_{10}
 +3\mu_{6}^{2})+2B_{00}\psi_{1}^{3})/(2\psi_{1}^{2}),
\end{array}\label{tran:11}
\end{equation}
and satisfying the condition $(\psi)_{yy}=\psi_{yy}$, is
involutive.


\subsubsection{Case $\mu_{5}=0$}

Equating the mixed derivatives $(\psi_{1xxx})_{y}=(\psi_{1y})_{xxx}$,
one gets
\begin{equation}
\begin{array}{rl}
B_{01y}= & (-12A_{00x}B_{2x}+12B_{02x}-4B_{10x}B_{2}-3B_{2x}A_{00}^{2}+8B_{2x}B_{10}\\
 & +3A_{00}B_{02}-A_{00}B_{10}B_{2})/18.
\end{array}\label{suff:20}
\end{equation}

The equation $(\psi_{0xxx})_{y}-(\psi_{0y})_{xxx}=0$ is a quadratic
algebraic equation with respect to $\psi_{1xx}$. Differentiating
this equation with respect to $y$ and $x$, and excluding from them
$\psi_{1xx}^{2}$ using it, one obtains the conditions
\begin{equation}
\begin{array}{rl}
36B_{02xx}= & (36A_{00xx}B_{2x}-18A_{00x}B_{2x}A_{00}+36A_{00x}B_{02}-12A_{00x}B_{10}B_{2}\\
 & +216B_{00yy}+72B_{00y}B_{2}-18B_{02x}A_{00}+12B_{10xx}B_{2}+24B_{10x}B_{2x}\\
 & +6B_{10x}A_{00}B_{2}-9B_{2x}A_{00}^{3}+36B_{2x}A_{00}B_{10}-72B_{2x}B_{01}\\
 & +72B_{2y}B_{00}+9A_{00}^{2}B_{02}-3A_{00}^{2}B_{10}B_{2}-30B_{02}B_{10}+10B_{10}^{2}B_{2}).
\end{array}\label{suff:21}
\end{equation}
\begin{equation}
\begin{array}{c}
2520\mu_{7}(2\psi_{1}\varphi_{x}^{2}\psi_{1xx}+9\varphi_{xx}^{2}\psi_{1}^{2}-20\varphi_{xx}\varphi_{x}\psi_{1x}\psi_{1}+8\varphi_{x}^{2}\psi_{1x}^{2})\\
+8\mu_{8}\psi_{1}\varphi_{x}(\psi_{1}\varphi_{xx}-\varphi_{x}\psi_{1x})+\mu_{9}\varphi_{x}^{2}\psi_{1}^{2}=0,
\end{array}\label{eq:11}
\end{equation}
 where
\[
\begin{array}{rl}
\mu_{7}= & -4A_{00xx}-6A_{00x}A_{00}+8B_{10x}-A_{00}^{3}+4A_{00}B_{10}-8B_{01},\\[2ex]
\mu_{8}= & 3456A_{00x}^{2}+1728A_{00x}A_{00}^{2}-8448A_{00x}B_{10}-38400B_{00y}+15360B_{01x}\\
 & -3840B_{10xx}-1920B_{10x}A_{00}-1320\mu_{7x}+216A_{00}^{4}-2112A_{00}^{2}B_{10}\\
 & +3840A_{00}B_{01}+855A_{00}\mu_{7}-12800B_{00}B_{2}+3456B_{10}^{2},\\
\mu_{9}= & 1503A_{00x}\mu_{7}-360\mu_{7xx}+1575\mu_{7x}A_{00}-\mu_{8x}-1278A_{00}^{2}\mu_{7}\\
 & +2A_{00}\mu_{8}-792B_{10}\mu_{7}.
\end{array}
\]

Considering equation (\ref{eq:11}), further study is separated in
two cases related with value of $\mu_{7}$, i.e., $\mu_{7}\neq0$
and $\mu_{7}=0$. 

\emph{3.2.2.1 Case $\mu_{7}\neq0$}

Substituting $\psi_{1xx}$, found from equation (\ref{eq:11}), into
$\psi_{1xxx}$, one gets the equation
\begin{equation}
\begin{array}{c}
8\varphi_{xx}\mu_{10}\psi_{1}-8\varphi_{x}\psi_{1x}\mu_{10}+\varphi_{x}\psi_{1}(-3333960A_{00x}A_{00}\mu_{7}^{2}-7560A_{00x}\mu_{7}\mu_{8}\\
+3810240B_{10x}\mu_{7}^{2}-2520\mu_{7x}\mu_{9}+2520\mu_{9x}\mu_{7}+952560A_{00}^{3}\mu_{7}^{2}\\
-793800A_{00}B_{10}\mu_{7}^{2}-2835A_{00}\mu_{7}\mu_{9}-2540160B_{01}\mu_{7}^{2}+2520B_{10}\mu_{7}\mu_{8}\\
+952560\mu_{7}^{3}+\mu_{8}\mu_{9})=0,
\end{array}\label{eq:12}
\end{equation}
 where
\[
\begin{array}{rl}
\mu_{10}= & -6191640A_{00x}\mu_{7}^{2}-2520\mu_{7x}\mu_{8}+2520\mu_{8x}\mu_{7}+1131165A_{00}^{2}\mu_{7}^{2}\\
 & -1890A_{00}\mu_{7}\mu_{8}+952560B_{10}\mu_{7}^{2}-630\mu_{7}\mu_{9}+\mu_{8}^{2}.
\end{array}
\]


\textit{Case $\mu_{10}\neq0$}

From equation (\ref{eq:12}), one arrives at
\begin{equation}
\varphi_{xx}=(\varphi_{x}\psi_{1x}+\varphi_{x}\mu_{11}\psi_{1})/\psi_{1},\label{tran:13}
\end{equation}
 where
\[
\begin{array}{rl}
8\mu_{10}\mu_{11}= & (3333960A_{00x}A_{00}\mu_{7}^{2}+7560A_{00x}\mu_{7}\mu_{8}-3810240B_{10x}\mu_{7}^{2}+2520\mu_{7x}\mu_{9}\\
 & -2520\mu_{9x}\mu_{7}-952560A_{00}^{3}\mu_{7}^{2}+793800A_{00}B_{10}\mu_{7}^{2}+2835A_{00}\mu_{7}\mu_{9}\\
 & +2540160B_{01}\mu_{7}^{2}-2520B_{10}\mu_{7}\mu_{8}-952560\mu_{7}^{3}-\mu_{8}\mu_{9}).
\end{array}
\]
 Substituting $\varphi_{xx}$ into $\varphi_{xxx}$, one gets the
condition
\begin{equation}
3\mu_{11y}+B_{2x}=0.\label{suff:23}
\end{equation}
 Differentiating equation (\ref{tran:13}) with respect to $y$, one
obtains the condition
\begin{equation}
\begin{array}{rl}
 & -7560A_{00x}\mu_{7}+7560A_{00}\mu_{11}\mu_{7}+2520B_{10}\mu_{7}+40320\mu_{11}^{2}\mu_{7}\\
 & +8\mu_{11}\mu_{8}+\mu_{9}-20160\mu_{7}\mu_{11x}=0.
\end{array}\label{suff:24}
\end{equation}

 The obtained system of equations for the functions $\varphi(x)$,
$\psi_{1}(x,y)$ and $\psi_{0}(x,y)$ consisting of the equations
\begin{equation}
\psi_{0y}=\psi_{1}(A_{00}+3\mu_{11}),\label{tran:14}
\end{equation}
\begin{equation}
\begin{array}{rl}
\psi_{0xxx}= & (20160\psi_{0xx}\psi_{1x}\mu_{7}\psi_{1}+20160\psi_{0xx}\mu_{11}\mu_{7}\psi_{1}^{2}-10080\psi_{0x}\psi_{1x}^{2}\mu_{7}\\
 & -20160\psi_{0x}\psi_{1x}\mu_{11}\mu_{7}\psi_{1}+\psi_{0x}\psi_{1}^{2}(-2520A_{00x}\mu_{7}\\
 & +2520A_{00}\mu_{11}\mu_{7}+840B_{10}\mu_{7}-30240\mu_{11}^{2}\mu_{7}-8\mu_{11}\mu_{8}\\
 & -\mu_{9})+6720B_{00}\mu_{7}\psi_{1}^{3})/(6720\mu_{7}\psi_{1}^{2}),
\end{array}\label{tran:15}
\end{equation}
\begin{equation}\label{tran:12}
    \begin{array}{rl}
    \psi_{1xx} =&( - 22680 \varphi_{xx}^2 \mu_{7} \psi_{1}^2 + 50400 \varphi_{xx} \varphi_{x} \psi_{1x} \mu_{7} \psi_{1} - 8 \varphi_{xx} \varphi_{x} \mu_{8} \psi_{1}^2\\
    &- 20160 \varphi_{x}^2 \psi_{1x}^2 \mu_{7} + 8 \varphi_{x}^2 \psi_{1x} \mu_{8} \psi_{1} - \varphi_{x}^2 \mu_{9} \psi_{1}^2)/(5040 \varphi_{x}^2 \mu_{7} \psi_{1}),
    \end{array}
\end{equation}
and satisfying the conditions $(\psi_{1xx})_{y}=(\psi_{1y})_{xx}$,
is involutive.


{\it Case $\mu_{10}=0$}

Substituting $\mu_{10}=0$ into equation (\ref{eq:12}), one obtains
the condition
\begin{equation}
\begin{array}{rl}
B_{10x}= & (3333960A_{00x}A_{00}\mu_{7}^{2}+7560A_{00x}\mu_{7}\mu_{8}+2520\mu_{7x}\mu_{9}-2520\mu_{9x}\mu_{7}\\
 & -952560A_{00}^{3}\mu_{7}^{2}+793800A_{00}B_{10}\mu_{7}^{2}+2835A_{00}\mu_{7}\mu_{9}+2540160B_{01}\mu_{7}^{2}\\
 & -2520B_{10}\mu_{7}\mu_{8}-952560\mu_{7}^{3}-\mu_{8}\mu_{9})/(3810240\mu_{7}^{2})
\end{array}\label{suff:22}
\end{equation}
 After all calculations the equations for the functions
$\varphi(x)$, $\psi_{1}(x,y)$ and $\psi_{0}(x,y)$ consisting of the equations
\begin{equation}
\psi_{0y}=(3\varphi_{xx}\psi_{1}-3\varphi_{x}\psi_{1x}+\varphi_{x}A_{00}\psi_{1})/\varphi_{x},\label{tran:16}
\end{equation}
\begin{equation}
\begin{array}{rl}
\psi_{0xxx}= & (-30240\varphi_{xx}^{2}\psi_{0x}\mu_{7}\psi_{1}^{2}+20160\varphi_{xx}\varphi_{x}\psi_{0xx}\mu_{7}\psi_{1}^{2}+40320\varphi_{xx}\varphi_{x}\psi_{0x}\psi_{1x}\mu_{7}\psi_{1}\\
 & +8\varphi_{xx}\varphi_{x}\psi_{0x}\psi_{1}^{2}(315A_{00}\mu_{7}-\mu_{8})-20160\varphi_{x}^{2}\psi_{0x}\psi_{1x}^{2}\mu_{7}\\
 & +8\varphi_{x}^{2}\psi_{0x}\psi_{1x}\psi_{1}(-315A_{00}\mu_{7}+\mu_{8})+\varphi_{x}^{2}\psi_{0x}\psi_{1}^{2}(-2520A_{00x}\mu_{7}\\
 & +840B_{10}\mu_{7}-\mu_{9})+6720\varphi_{x}^{2}B_{00}\mu_{7}\psi_{1}^{3})/(6720\varphi_{x}^{2}\mu_{7}\psi_{1}^{2}),
\end{array}\label{tran:17}
\end{equation}
\begin{equation}
\begin{array}{rl}
\varphi_{xxx}= & (-10080\varphi_{xx}^{2}\mu_{7}\psi_{1}^{2}+40320\varphi_{xx}\varphi_{x}\psi_{1x}\mu_{7}\psi_{1}+8\varphi_{xx}\varphi_{x}\psi_{1}^{2}(315A_{00}\mu_{7}-\mu_{8})\\
 & -20160\varphi_{x}^{2}\psi_{1x}^{2}\mu_{7}+8\varphi_{x}^{2}\psi_{1x}\psi_{1}(-315A_{00}\mu_{7}+\mu_{8})+\varphi_{x}^{2}\psi_{1}^{2}(-2520A_{00x}\mu_{7}\\
 & +840B_{10}\mu_{7}-\mu_{9}))/(6720\varphi_{x}\mu_{7}\psi_{1}^{2}),
\end{array}\label{tran:18}
\end{equation}
\begin{equation}\label{tran:27}
    \begin{array}{rl}
    \psi_{1xx} =&( - 22680 \varphi_{xx}^2 \mu_{7} \psi_{1}^2 + 50400 \varphi_{xx} \varphi_{x} \psi_{1x} \mu_{7} \psi_{1} - 8 \varphi_{xx} \varphi_{x} \mu_{8} \psi_{1}^2\\
    &- 20160 \varphi_{x}^2 \psi_{1x}^2 \mu_{7} + 8 \varphi_{x}^2 \psi_{1x} \mu_{8} \psi_{1} - \varphi_{x}^2 \mu_{9} \psi_{1}^2)/(5040 \varphi_{x}^2 \mu_{7} \psi_{1}),
    \end{array}
\end{equation}
and satisfying the conditions $(\psi_{1xx})_{y}=(\psi_{1y})_{xx}$,
is involutive.


\emph{3.2.2.2 Case $\mu_{7}=0$}

Substituting $\mu_{7}=0$ into equation (\ref{eq:11}), one gets
\begin{equation}
8\varphi_{xx}\mu_{8}\psi_{1}-8\varphi_{x}\psi_{1x}\mu_{8}+\varphi_{x}\psi_{1}(-\mu_{8x}+2A_{00}\mu_{8})=0.\label{eq:13}
\end{equation}
 Further analysis of the compatibility depends upon the value of
$\mu_{8}$.

{\it Case $\mu_{8}\neq0$}

Solving equation (\ref{eq:13}) with respect to $\varphi_{xx}$, and substituting it into $\varphi_{xxx}$, one obtains
\begin{equation}
\begin{array}{rl}
\psi_{1xx}= & (192\psi_{1x}^{2}\mu_{8}^{2}+16\psi_{1x}\mu_{8}\psi_{1}(\mu_{8x}-2A_{00}\mu_{8})+\psi_{1}^{2}(-64A_{00x}\mu_{8}^{2}-64\mu_{8xx}\mu_{8}\\
 & +71\mu_{8x}^{2}-4\mu_{8x}A_{00}\mu_{8}-20A_{00}^{2}\mu_{8}^{2}+64B_{10}\mu_{8}^{2}))/(128\mu_{8}^{2}\psi_{1}),
\end{array}\label{tran:20}
\end{equation}
Substituting $\psi_{1xx}$ into the found earlier expression for $\psi_{1xxx}$, one arrives at the condition
\begin{equation}
\begin{array}{rl}
\mu_{8xxx}= & (48A_{00x}\mu_{8x}\mu_{8}^{2}+96A_{00x}A_{00}\mu_{8}^{3}-128B_{10x}\mu_{8}^{3}+300\mu_{8xx}\mu_{8x}\mu_{8}\\
 & -225\mu_{8x}^{3}+12\mu_{8x}A_{00}^{2}\mu_{8}^{2}-32\mu_{8x}B_{10}\mu_{8}^{2}+24A_{00}^{3}\mu_{8}^{3}\\
 & -96A_{00}B_{10}\mu_{8}^{3}+192B_{01}\mu_{8}^{3})/(80\mu_{8}^{2}),
\end{array}\label{suff:25}
\end{equation}

The obtained system of equations for the functions $\varphi(x)$,
$\psi_{1}(x,y)$ and $\psi_{0}(x,y)$ consisting of the equations
\begin{equation}
\psi_{0y}=(\psi_{1}(3\mu_{8x}+2A_{00}\mu_{8}))/(8\mu_{8}),\label{tran:21}
\end{equation}
\begin{equation}
\begin{array}{rl}
\psi_{0xxx}= & (384\psi_{0xx}\psi_{1x}\mu_{8}^{2}\psi_{1}+48\psi_{0xx}\mu_{8}\psi_{1}^{2}(\mu_{8x}-2A_{00}\mu_{8})-192\psi_{0x}\psi_{1x}^{2}\mu_{8}^{2}\\
 & +48\psi_{0x}\psi_{1x}\mu_{8}\psi_{1}(-\mu_{8x}+2A_{00}\mu_{8})+\psi_{0x}\psi_{1}^{2}(-96A_{00x}\mu_{8}^{2}-48\mu_{8xx}\mu_{8}\\
 & +51\mu_{8x}^{2}+12\mu_{8x}A_{00}\mu_{8}-36A_{00}^{2}\mu_{8}^{2}
 +64B_{10}\mu_{8}^{2})\\
 &+128B_{00}\mu_{8}^{2}\psi_{1}^{3})/(128\mu_{8}^{2}\psi_{1}^{2}),
\end{array}\label{tran:22}
\end{equation}
\begin{equation}\label{tran:19}
            \varphi_{xx} =  (8 \varphi_{x} \psi_{1x} \mu_{8} + \varphi_{x} \psi_{1} (\mu_{8x} - 2 A_{00} \mu_{8}))/(8 \mu_{8}
            \psi_{1})
        \end{equation}
is
involutive.


{\it Case $\mu_{8}=0$}

 Since $\mu_{8}=0$, then equation (\ref{eq:13}) is satisfied. The obtained system of equations for the functions $\varphi(x)$, $\psi_{1}(x,y)$ and $\psi_{0}(x,y)$ consisting of the equations
\begin{equation}
\psi_{0y}=(3\varphi_{xx}\psi_{1}-3\varphi_{x}\psi_{1x}+\varphi_{x}A_{00}\psi_{1})/\varphi_{x},\label{tran:23}
\end{equation}
\begin{equation}
\begin{array}{rl}
\psi_{0xxx}= & (-9\varphi_{xx}^{2}\psi_{0x}\psi_{1}+24\varphi_{xx}\varphi_{x}\psi_{0xx}\psi_{1}-12\varphi_{xx}\varphi_{x}\psi_{0x}\psi_{1x}\\
 & +3\varphi_{xx}\varphi_{x}\psi_{0x}A_{00}\psi_{1}+6\varphi_{x}^{2}\psi_{0x}\psi_{1xx}-3\varphi_{x}^{2}\psi_{0x}\psi_{1x}A_{00}\\
 & +\varphi_{x}^{2}\psi_{0x}\psi_{1}(-3A_{00x}+B_{10})+8\varphi_{x}^{2}B_{00}\psi_{1}^{2})/(8\varphi_{x}^{2}\psi_{1}),
\end{array}\label{tran:24}
\end{equation}
\begin{equation}
\begin{array}{rl}
\psi_{1xxx}= & (90\varphi_{xx}^{3}\psi_{1}^{2}-540\varphi_{xx}^{2}\varphi_{x}\psi_{1x}\psi_{1}+135\varphi_{xx}^{2}\varphi_{x}A_{00}\psi_{1}^{2}+420\varphi_{xx}\varphi_{x}^{2}\psi_{1xx}\psi_{1}\\
 & +240\varphi_{xx}\varphi_{x}^{2}\psi_{1x}^{2}-330\varphi_{xx}\varphi_{x}^{2}\psi_{1x}A_{00}\psi_{1}+3\varphi_{xx}\varphi_{x}^{2}\psi_{1}^{2}(-14A_{00x}+19A_{00}^{2}\\
 & -14B_{10})-120\varphi_{x}^{3}\psi_{1xx}\psi_{1x}+90\varphi_{x}^{3}\psi_{1xx}A_{00}\psi_{1}+60\varphi_{x}^{3}\psi_{1x}^{2}A_{00}\\
 & +\varphi_{x}^{3}\psi_{1x}\psi_{1}(12A_{00x}-57A_{00}^{2}+52B_{10})+\varphi_{x}^{3}\psi_{1}^{2}(-21A_{00x}A_{00}\\
 & +24B_{10x}+6A_{00}^{3}-5A_{00}B_{10}-16B_{01}))/(80\varphi_{x}^{3}\psi_{1}),
\end{array}\label{tran:25}
\end{equation}
\begin{equation}
\begin{array}{rl}
\varphi_{xxx}= & (15\varphi_{xx}^{2}\psi_{1}-12\varphi_{xx}\varphi_{x}\psi_{1x}+3\varphi_{xx}\varphi_{x}A_{00}\psi_{1}+6\varphi_{x}^{2}\psi_{1xx}\\
 & -3\varphi_{x}^{2}\psi_{1x}A_{00}+\varphi_{x}^{2}\psi_{1}(-3A_{00x}+B_{10}))/(8\varphi_{x}\psi_{1})
\end{array}\label{tran:26}
\end{equation}
is involutive.

\section{Linearization theorems}

All obtained results can be summarized in the following theorems.
\begin{thm}
\label{thm1}

Any fourth-order ordinary differential equation obtained from a linear
equation (\ref{eq:03}) by a tangent transformation (\ref{eq:02})
has to be either to the form (\ref{eq:08}) or (\ref{eq:09}).
\end{thm}

\begin{thm}
\label{thm2}

Sufficient conditions for equation (\ref{eq:08}) to be linearizable
via a tangent transformation with $\varphi_{p}=\varphi_{y}=0$ and
$\alpha=0$ are as follows.
\begin{description}
\item [{{{{(a)}}}}] If $A_{1}\neq0$, then the conditions are
(\ref{suff:01}), (\ref{suff:02}), (\ref{suff:03}),
(\ref{suff:04}), (\ref{suff:06}), (\ref{suff:07}) and
(\ref{suff:08}).
\item [{{{{(b)}}}}] If $A_{1}=0$ and $\mu_{5}\neq0$, then the conditions
are (\ref{suff:09}), (\ref{suff:11}), (\ref{suff:10}),
(\ref{suff:12}), (\ref{suff:13}), (\ref{suff:14}),
(\ref{suff:15}), (\ref{suff:16}), (\ref{suff:17}), (\ref{suff:26}) and
(\ref{suff:19}).
\item [{{{{(c)}}}}] If $A_{1}=0,\mu_{5}=0,\mu_{7}\neq0$ and $\mu_{10}\neq0$,
then the conditions are (\ref{suff:09}), (\ref{suff:11}),
(\ref{suff:10}), (\ref{suff:12}), (\ref{suff:13}),
(\ref{suff:14}), (\ref{suff:15}), (\ref{suff:20}),
(\ref{suff:21}), (\ref{suff:23}) and (\ref{suff:24}).
\item [{{{{(d)}}}}] If $A_{1}=0,\mu_{5}=0,\mu_{7}\neq0$ and $\mu_{10}=0$,
then the conditions are (\ref{suff:09}), (\ref{suff:11}),
(\ref{suff:10}), (\ref{suff:12}), (\ref{suff:13}),
(\ref{suff:14}), (\ref{suff:15}), (\ref{suff:20}), (\ref{suff:21})
and (\ref{suff:22}).
\item [{{{{(e)}}}}] If $A_{1}=0,\mu_{5}=0,\mu_{7}=0$ and $\mu_{8}\neq0$,
then the conditions are (\ref{suff:09}), (\ref{suff:11}),
(\ref{suff:10}), (\ref{suff:12}), (\ref{suff:13}),
(\ref{suff:14}), (\ref{suff:15}), (\ref{suff:20}),
(\ref{suff:21}),  and (\ref{suff:25}).
\item [{{{{(f)}}}}] If $A_{1}=0,\mu_{5}=0,\mu_{7}=0$ and $\mu_{8}=0$,
then the conditions are (\ref{suff:09}), (\ref{suff:11}),
(\ref{suff:10}), (\ref{suff:12}), (\ref{suff:13}),
(\ref{suff:14}), (\ref{suff:15}), (\ref{suff:20}) and
(\ref{suff:21}).
\end{description}
\end{thm}

\begin{thm}
\label{thm3}

Provided that the sufficient conditions in Theorem \ref{thm2} are
satisfied, the transformation (\ref{eq:02}) mapping equation (\ref{eq:08})
to a linear equation $u'''=0$ is obtained by solving the following
compatible system of equations for the functions $\varphi(x)$ and
$\psi(x,y,p)$:
\begin{description}
\item [{{{{(a)}}}}] (\ref{tran:02}), (\ref{tran:03}), (\ref{tran:04}), (\ref{tran:05})
and (\ref{tran:01}),
\item [{{{{(b)}}}}]  (\ref{tran:06}), (\ref{tran:07}),
(\ref{tran:08}), (\ref{tran:09}), (\ref{tran:10}) and (\ref{tran:11}),
\item [{{{{(c)}}}}] (\ref{tran:06}), (\ref{tran:07}), (\ref{tran:13}),
(\ref{tran:14}), (\ref{tran:15}) and (\ref{tran:12}),
\item [{{{{(d)}}}}] (\ref{tran:06}), (\ref{tran:07}), (\ref{tran:16}), (\ref{tran:17}), (\ref{tran:18}) and
(\ref{tran:27}),
\item [{{{{(e)}}}}] (\ref{tran:06}), (\ref{tran:07}),
(\ref{tran:20}), (\ref{tran:21}), (\ref{tran:22}) and
(\ref{tran:19}),
\item [{{{{(f)}}}}] (\ref{tran:06}), (\ref{tran:07}), (\ref{tran:23}),
(\ref{tran:24}), (\ref{tran:25}) and (\ref{tran:26}).
\end{description}
\end{thm}

\section{Examples}

{\bf Example 1}. Consider the nonlinear ordinary differential
equation
\begin{equation}
\label{ex:01}
xy{y^{\left( 4 \right)}} + \left( {4xy' + 3y} \right)y''' + 3x{y''^2} + 9y'y'' = 0.
\end{equation}
It is an equation of the form (\ref{eq:08}) with the coefficients
\begin{equation*}
\begin{array}{cc}
&A_1 =  0, \, ~  A_0 =\frac{{{\rm{4px  +  3y}}}}{{xy}}, \, ~  B_3 = 0, \, ~
 \\&
 B_2 = \frac{3}{y}, \, ~B_1 = \frac{{9p}}{{xy}},\, ~
 B_0 =0,\,\, ~  q = 0,\,\, ~\lambda=0.
\end{array}
\end{equation*}
The other functions corresponding to the case where $A_{1}=0$ are
\[
  \mu_1 = \frac{1}{y}, \, ~ \mu_2 = 0, \, ~ A_{00} =  \frac{3}{x}, \, ~  B_{10} = 0, \, ~ B_{00} = 0, \, ~B_{01} = 0,\, ~
 B_{02} =0,
\]
\[
    \mu_5 =0, \, ~ \mu_7 = \frac{3}{x^3}, \, ~\mu_8 =  \frac{21519}{x^4}, \, ~  \mu_9 = \frac{111672}{x^5},
\]
\[
    \mu_{10} = -\frac{17915904}{x^8}, \, ~\mu_{11} = -\frac{1}{x}.
\]
Since $\mu_{5}=0$, $\mu_{7}\neq 0$ and $\mu_{10}\neq 0$, and all these coefficients obey the
conditions
 (\ref{suff:09}), (\ref{suff:11}),
(\ref{suff:10}), (\ref{suff:12}), (\ref{suff:13}),
(\ref{suff:14}), (\ref{suff:15}), (\ref{suff:20}),
(\ref{suff:21}), (\ref{suff:23}) and (\ref{suff:24}), one concludes that
equation (\ref {ex:01}) is linearizable.
Applying Theorem \ref{thm3}, the linearizing transformation is found by solving the following equations.
For the function $\psi_{1}(x,y)$ one has
\begin{equation}
\label{ex:03}
\psi_{1y} = \frac{\psi_{1}}{y}, \qquad \psi_{1xx} = \frac{3 \psi_{1x}^2 x^2 - 2 \psi_{1x} \psi_{1} x - \psi_{1}^2}{2 \psi_{1} x^2}.
\end{equation}
The function $\psi_{1}= x y$ is a particular solution of
equations (\ref {ex:03}). The equations for the function  $\psi_{0}(x,y)$ become
\[
\psi_{0y}=0, \qquad \psi_{0xxx}=0.
\]
One can choose the simplest solution of these equations: $\psi_{0}=0$.
For finding the function $\varphi(x)$ one has to solve the equation
\begin{equation}
\label{ex:02}
\varphi_{xx}=0.
\end{equation}
Choosing the particular solution $\varphi=x$, one obtains the
linearizing transformation
\begin{equation}
\label{ex:04}
t = x, \qquad u = x y p.
\end{equation}
Thus, the nonlinear equation (\ref{ex:01}) can be mapped by transformation
(\ref{ex:04}) into the linear equation $u^{(3)}=0$.


\noindent {\bf Example 2}. Consider the nonlinear ordinary differential
equation
\begin{equation}
\label{ex:05}
{y^3}{y^{\left( 4 \right)}} - 4{y^2}y'y''' - 3{y^2}{y''^2} + 12y{y'^2}y'' - 6{y'^4} = 0.
\end{equation}
It is an equation of form (\ref{eq:08}) with the coefficients
\begin{equation*}
\begin{array}{cc}
&A_1 =  0, \, ~  A_0 =-\frac{{{\rm{4p}}}}{{y}}, \, ~  B_3 = 0, \, ~
 \\&
 B_2 = -\frac{3}{y}, \, ~B_1 = \frac{{12p^2}}{{y^2}},\, ~
 B_0 =-\frac{6p^4}{y^3},\,\, ~  q = 0,\,\, ~\lambda=0.
\end{array}
\end{equation*}
The other coefficients applying in the case $A_1 = 0$ are
\begin{equation*}
    \mu_1 =- \frac{1}{y}, \, ~ \mu_2 = 0, \, ~ A_{00} = 0, \, ~  B_{10} = 0, \, ~ B_{00} = 0,
 \end{equation*}
\begin{equation*}
   B_{01} = 0, \, ~B_{02} =0,  \, ~\mu_5 =0, \, ~ \mu_7 =0, \, ~\mu_8 = 0, \, ~  \mu_9 = 0.
\end{equation*}
Since $\mu_{5}=0$, $\mu_{7}=0$ and $\mu_{8}= 0$, one can check that these coefficients obey the conditions
(\ref{suff:09}), (\ref{suff:11}),
(\ref{suff:10}), (\ref{suff:12}), (\ref{suff:13}),
(\ref{suff:14}), (\ref{suff:15}), (\ref{suff:20}) and
(\ref{suff:21}). Thus,
equation (\ref {ex:05}) is linearizable.
The linearizing transformation is found as follows. According to Theorem
\ref{thm3}, one obtains that the function $\psi_{1}$ satisfies the equation
\begin{equation*}
\psi_{1y}=-\frac{\psi_{1}}{y}.
\end{equation*}
One can take the  solution $\psi_{1}=\frac{1}{y}$ of
this equation. For the function $\varphi(x)$ the equations are
\begin{equation}
\label{ex:06}
\varphi_{xxx} = \frac{15\varphi_{xx}^2}{8\varphi_x}, \qquad \varphi_{xx}^3 =0.
\end{equation}
The function $\varphi = x$ is a particular solution of equations (\ref {ex:06}).
For finding the function $\psi_{0}(x,y)$, one has to solve the  system of equations
\begin{equation}
\label{ex:07}
\psi_{0y}=0, \qquad \psi_{0xxx}=0.
\end{equation}
A particular solution of these equations is $\psi_{0}=0$. Hence one obtains the
linearizing transformation
\begin{equation}
\label{ex:08}
t = x, \qquad u = p/y.
\end{equation}
which maps nonlinear equation (\ref{ex:05}) into the linear equation $u^{(3)}=0$.


%


\section{Acknowledgement}

This research was financially supported by the Discovery Based Development
Grant (P-10-11295), National Science and Technology Development Agency,
Thailand and Naresuan University. The research of SVM was partially supported by the Office of the Higher Education
Commission under NRU project (SUT).

\appendix

\section{Appendix}

\label{Appendix}

For proving theorems one needs relations between $\varphi(x,y,p)$
and $\psi(x,y,p)$ and coefficients of equation (\ref{eq:08}).
These relations are presented here.
\begin{align}
A_{1} = &((3 \psi_{pp} p - 4 \psi_{pp}) \varphi_{y} + 3
\varphi_{x} \psi_{pp})/((\varphi_{x} + \varphi_{y} p)
\psi_{pp}),\label{ap:01}
\end{align}
\begin{align}
A_{0} = &( - (3 \varphi_{xx} + \varphi_{yy} p^2 + 2 \varphi_{xy}
p)
\psi_{pp} - (3 \psi_{py} p + \psi_{y} + 3 \psi_{px}) \varphi_{x} \nonumber\\
&- (3 \psi_{py} p^2 - \psi_{x} + 3 \psi_{px} p)
\varphi_{y}))/((\varphi_{x} + \varphi_{y} p)
\psi_{pp}),\label{ap:02}
\end{align}
\begin{align}
B_{3} = &(((2 \psi_{ppp} p - 3 \psi_{pp}) \varphi_{y} +
\varphi_{x}
\psi_{ppp}) \varphi_{x} - (3 (\psi_{pp} p - \psi_{pp}) \nonumber\\
&- \psi_{ppp} p^2) \varphi_{y}^2)/((\varphi_{x} + \varphi_{y} p)^2
\psi_{pp}),\label{ap:03}
\end{align}
\begin{align}
B_{2} = &( - 3 (((\psi_{pp} p - 2 \psi_{pp}) \varphi_{y} +
\varphi_{x} \psi_{pp}) \varphi_{xx} + (\psi_{pp} p - \psi_{pp})
\varphi_{yy} \varphi_{y} p^2 \nonumber\\
&+ ((2 \psi_{pp} p + \psi_{pp}) \varphi_{x} + (2 \psi_{pp} p - 3
\psi_{pp}) \varphi_{y} p) \varphi_{xy} - (\psi_{ppx} + \psi_{ppy}
p
+ \psi_{py}) \varphi_{x}^2 \nonumber\\
&- ((2 \psi_{ppy} p^2 - \psi_{y} + 2 \psi_{ppx} p - 2 \psi_{px})
\varphi_{y} - (\psi_{pp} p + \psi_{pp}) \varphi_{yy} p)
\varphi_{x}
\nonumber\\
&- (\psi_{ppy} p^3 + \psi_{x} + \psi_{ppx} p^2 - \psi_{py} p^2 - 2
\psi_{px} p) \varphi_{y}^2))/((\varphi_{x} + \varphi_{y} p)^2
\psi_{pp}),\label{ap:04}
\end{align}
\begin{align}
B_{1} = &( - ((\varphi_{yyy} \varphi_{y} - 3 \varphi_{yy}^2) p^4 -
3 \varphi_{xx}^2 - 12 \varphi_{xy}^2 p^2 + (\varphi_{xxx} + 3
\varphi_{xxy} p \nonumber\\
&+ 3 \varphi_{xyy} p^2) (\varphi_{x} + \varphi_{y} p)) \psi_{pp} +
3 (2 \psi_{py} p^2 - \psi_{x} + 2 \psi_{px} p) \varphi_{yy}
\varphi_{y} p^2 \nonumber\\
&- 3 (\psi_{xy} + \psi_{yy} p + \psi_{pyy} p^2 + \psi_{pxx} + 2
\psi_{pxy} p) \varphi_{x}^2 + 3 (\psi_{xy} p + \psi_{xx} -
\psi_{pyy} p^3 \nonumber\\
&- \psi_{pxx} p - 2 \psi_{pxy} p^2) \varphi_{y}^2 p + 3 ((2
\psi_{py} p + \psi_{y} + 2 \psi_{px}) \varphi_{x} - 2 \varphi_{yy}
\psi_{pp} p^2 \nonumber\\
&- (2 \psi_{x} + \psi_{y} p - 2 \psi_{py} p^2 - 2 \psi_{px} p)
\varphi_{y}) \varphi_{xx} + 3 ((\psi_{x} + 3 \psi_{y} p + 4
\psi_{py} p^2 \nonumber\\
&+ 4 \psi_{px} p) \varphi_{x} - 4 \varphi_{xx} + \varphi_{yy} p^2)
\psi_{pp} p - (3 \psi_{x} + \psi_{y} p - 4 \psi_{py} p^2 - 4
\psi_{px} p) \varphi_{y} p) \varphi_{xy} \nonumber\\
&+ (3 (\psi_{xx} - \psi_{yy} p^2 - 2 \psi_{pyy} p^3 - 2 \psi_{pxx}
p
- 4 \psi_{pxy} p^2) \varphi_{y} + (3 (\psi_{x} + 2 \psi_{y} p \nonumber\\
&+ 2 \psi_{py} p^2 + 2 \psi_{px} p) \varphi_{yy} + \varphi_{yyy}
\psi_{pp} p^2) p) \varphi_{x}))/((\varphi_{x} + \varphi_{y} p)^2
\psi_{pp}),\label{ap:05}
\end{align}
\begin{align}
B_{0} = &( - ((\varphi_{yyy} \varphi_{y} - 3 \varphi_{yy}^2) p^4 -
3 \varphi_{xx}^2 - 12 \varphi_{xy}^2 p^2 + (\varphi_{xxx} + 3
\varphi_{xxy} p \nonumber\\
&+ 3 \varphi_{xyy} p^2) (\varphi_{x} + \varphi_{y} p)) (\psi_{x} +
\psi_{y} p) + 3 (\psi_{xx} + \psi_{yy} p^2 + 2 \psi_{xy} p)
\varphi_{yy} \varphi_{y} p^3 \nonumber\\
&+ 3 ((\psi_{xx} + \psi_{yy} p^2 + 2 \psi_{xy} p) (\varphi_{x} +
\varphi_{y} p) - 2 (\psi_{x} + \psi_{y} p) \varphi_{yy} p^2)
\varphi_{xx} \nonumber\\
&+ 6 ((\psi_{xx} + \psi_{yy} p^2 + 2 \psi_{xy} p) (\varphi_{x} +
\varphi_{y} p) - 2 \varphi_{xx} + \varphi_{yy} p^2) (\psi_{x} +
\psi_{y} p)) \varphi_{xy} p \nonumber\\
&- ((3 \psi_{xxy} + \psi_{yyy} p^2) p + \psi_{xxx} + 3 \psi_{xyy}
p^2) (\varphi_{x}^2 + \varphi_{y}^2 p^2) \nonumber\\
&- (2 ((3 \psi_{xxy} + \psi_{yyy} p^2) p + \psi_{xxx} + 3
\psi_{xyy}
p^2) \varphi_{y} - (3 (\psi_{xx} + \psi_{yy} p^2 \nonumber\\
&+ 2 \psi_{xy} p) \varphi_{yy} + (\psi_{x} + \psi_{y} p)
\varphi_{yyy} p) p) \varphi_{x} p))/((\varphi_{x} + \varphi_{y}
p)^2 \psi_{pp}).\label{ap:06}
\end{align}
\begin{equation}
\begin{array}{rl}
A_{0xy}= & (-18A_{0x}A_{1}\mu_{1}p^{2}-6A_{0y}A_{0}A_{1}p^{2}+18A_{0y}A_{1}\mu_{1}p^{3}+6A_{0y}\mu_{2}p^{2}\\
 & +12B_{0p}A_{1}p+6B_{1x}A_{1}p-6B_{1y}A_{1}p^{2}+36\mu_{1xy}A_{1}p^{3}\\
 & -18\mu_{1x}A_{0}A_{1}p^{2}+72\mu_{1x}A_{1}\mu_{1}p^{3}+6\mu_{1yy}A_{1}p^{4}+12\mu_{1y}A_{0}A_{1}p^{3}
 \\
 & -18\mu_{1y}A_{1}\mu_{1}p^{4}-12\mu_{1y}\mu_{2}p^{3}+6\mu_{3x}A_{1}p-15\mu_{4p}A_{1}p
 \\
 & -6A_{0}^{2}A_{1}\mu_{1}p^{2}+2A_{0}A_{1}B_{1}p+30A_{0}A_{1}\mu_{1}^{2}p^{3}+8A_{0}A_{1}\mu_{3}p\\
 & +6A_{0}\mu_{1}\mu_{2}p^{2}+4A_{1}^{2}B_{0}p-5A_{1}^{2}\mu_{4}p-18A_{1}B_{0}-12A_{1}B_{1}\mu_{1}p^{2}\\
 & -30A_{1}\mu_{1}^{3}p^{4}-36A_{1}\mu_{1}\mu_{3}p^{2}+9A_{1}\mu_{4}-2B_{1}\mu_{2}p-12\mu_{1}^{2}\mu_{2}p^{3}
 \\
 & +4\mu_{2}\mu_{3}p)/(18A_{1}p^{2}),
\end{array}\label{suff:07}
\end{equation}
\begin{equation}
\begin{array}{rl}
\mu_{1yyy}= & (12A_{0x}A_{1}\mu_{3}+24A_{0y}A_{1}\mu_{3}p-2A_{1y}A_{1}B_{0}p+A_{1y}A_{1}\mu_{4}p-6B_{0px}A_{1}\\
 & -12B_{0py}A_{1}p-2B_{0p}A_{0}A_{1}-6B_{0p}A_{1}\mu_{1}p+2B_{0p}\mu_{2}-2B_{0x}A_{1}^{2}\\
 & -4B_{0y}A_{1}^{2}p+18B_{0y}A_{1}+6B_{1xy}A_{1}p+6B_{1x}A_{1}\mu_{1}p+6B_{1yy}A_{1}p^{2}\\
 & +2B_{1y}A_{0}A_{1}p+6B_{1y}A_{1}\mu_{1}p^{2}-2B_{1y}\mu_{2}p-12\mu_{1xyy}A_{1}p^{3}-36\mu_{1xy}A_{1}\mu_{1}p^{3}\\
 & -36\mu_{1x}\mu_{1y}A_{1}p^{3}+6\mu_{1x}A_{1}B_{1}p-36\mu_{1x}A_{1}\mu_{1}^{2}p^{3}-36\mu_{1x}A_{1}\mu_{3}p\\
 & -4\mu_{1yy}A_{0}A_{1}p^{3}-12\mu_{1yy}A_{1}\mu_{1}p^{4}+4\mu_{1yy}\mu_{2}p^{3}-18\mu_{1y}^{2}A_{1}p^{4}\\
 & -12\mu_{1y}A_{0}A_{1}\mu_{1}p^{3}+6\mu_{1y}A_{1}B_{1}p^{2}-24\mu_{1y}A_{1}\mu_{3}p^{2}+12\mu_{1y}\mu_{1}\mu_{2}p^{3}\\
 & -6\mu_{3x}\mu_{2}+3\mu_{4px}A_{1}+6\mu_{4py}A_{1}p+\mu_{4p}A_{0}A_{1}+3\mu_{4p}A_{1}\mu_{1}p-\mu_{4p}\mu_{2}\\
 & +\mu_{4x}A_{1}^{2}+2\mu_{4y}A_{1}^{2}p+4A_{0}^{2}A_{1}\mu_{3}+2A_{0}A_{1}B_{1}\mu_{1}p-4A_{0}A_{1}\mu_{1}^{3}p^{3}\\
 & -4A_{0}\mu_{2}\mu_{3}-2A_{1}^{2}B_{0}\mu_{1}p+A_{1}^{2}\mu_{1}\mu_{4}p-2A_{1}B_{0}B_{2}+24A_{1}B_{0}\mu_{1}\\
 & -12A_{1}B_{1}\mu_{3}+A_{1}B_{2}\mu_{4}+6A_{1}\mu_{1}^{4}p^{4}+12A_{1}\mu_{1}^{2}\mu_{3}p^{2}-3A_{1}\mu_{1}\mu_{4}\\
 & +24A_{1}\mu_{3}^{2}-2B_{1}\mu_{1}\mu_{2}p+4\mu_{1}^{3}\mu_{2}p^{3}+12\mu_{1}\mu_{2}\mu_{3}p)/(6A_{1}p^{4}).
\end{array}\label{suff:08}
\end{equation}


\end{document}